\documentclass[10pt,conference]{ieeeconf}
\IEEEoverridecommandlockouts
\usepackage{amsmath,amssymb,latexsym,epsfig,psfrag,graphicx}
\newtheorem{prp}{Proposition}
\newtheorem{thm}[prp]{Theorem}

\newenvironment{pf}{\noindent{\it Proof. }}{\hfill$\Box$\smallbreak}
\newenvironment{pf*}[1]{\smallbreak\noindent{\it #1}}{\hfill$\Box$\smallbreak}
\newcounter{definition}

\newcounter{remark}
\newenvironment{rmk}{\addtocounter{remark}{1}\smallbreak\noindent
  {\em Remark \theremark.}}{\smallbreak}
\newcommand{\realR}{\mathbb{R}}

\newcommand{\averageE}{\mathbb{E}}

\DeclareMathOperator{\argmax}{argmax}

\begin{document}
\title{On Minimax Optimal Dual Control for Fully Actuated Systems}
  \author{Anders Rantzer
  \thanks{The author is affiliated with Automatic Control LTH, Lund
    University, Box 118, SE-221 00 Lund, Sweden. He is a member of the Excellence Center ELLIIT and Wallenberg AI, Autonomous Systems and Software Program (WASP). Support was received from the European Research Council (Advanced Grant 834142) }}
\maketitle

\begin{abstract}%
A multi-variable adaptive controller is derived as the explicit solution to a minimax dynamic game. The minimizing player selects the control action as a function of past state measurements and inputs. The maximizing player selects disturbances and model parameters for the underlying linear time-invariant dynamics. This leads to a Bellman equation that can be solved explicitly for the case with unitary B-matrix known up to a sign and no input penalty. The minimizing policy is a dual controller that optimizes the tradeoff between exploration and exploitation.
\end{abstract}

\section{The Problem}
\label{sec:minimax}
This is a paper about adaptive control with worst-case models for disturbances and uncertain parameters, as discussed by \cite{cusumano1988nonlinear}, \cite{didinsky1994minimax} and \cite{sun1987theory}: 

\smallskip

\textbf{Problem:}
\emph{Given a number $\alpha>0$ and state dimension $n$, determine
\begin{align}
  J^*(x_0)&:=\inf_\mu\underbrace{\sup_{A,B}\sup_{w}\averageE\sum_{t=0}^\infty\left(|x_t|^2-\gamma^2|w_t|^2\right)}_{J_\mu(x_0)}
\label{eqn:infsup}
\end{align}
when $AA^\top=\alpha^2I$, $B=\pm I$, the state $x_t\in\realR^n$ and input $u_t\in\realR^n$ are determined from $w$ as
\begin{align}
  x_{t+1}&=Ax_t+Bu_t+w_t 
\label{eqn:plant}\\
  u_t&=\mu_t(x_0,\ldots,x_t,u_0,\ldots,u_{t-1}).\label{eqn:mu_LQ} 
\end{align}
and the control policy $\mu_t$ generates a random input $u_t$ based on measurements of $x_0,\ldots,x_t$
}

\smallskip

The problem is a dynamic game, where the $\mu$-player tries to minimize the cost, while the $(A,B,w)$-player tries to maximize it. Without uncertainty in $A$ and $B$, this would be the standard game formulation of $H_\infty$ optimal control \cite{Basar/B95} and $\gamma$ would serve as an upper bound on the gain from $w$ to $x$. In our formulation, the maximizing player can choose not only $w$, but also the $A$ and $B$. These parameters are unknown but constant, so an optimal feedback law tends to ``learn'' $A,B$ early on, in order to exploit this knowledge later. 
Such nonlinear adaptive controllers can stabilize and optimize the behavior also when no linear time-invariant controller can simultaneously stabilize for all parameter values. The parameter $\gamma$ still serves as a bound on the gain from disturbances to (average) state errors.

The main contribution of this paper is to formulate an exact solution to the Bellman equation that proves and generalizes a conjecture stated by Vinnicombe for scalar systems in \cite{vinnicombe2004examples}. Related work with suboptimal bounds has previously been published in \cite{rantzer2020ifac,orlov2018adaptive,rantzer2021minimax,kjellqvist2022learning,kjellqvist2022minimax}. 

\section{Notation}
The set of $n\times m$ matrices with real coefficients is denoted $\realR^{n\times m}$. The transpose of a matrix $A$ is denoted $A^\top$. For a symmetric matrix $A\in\realR^{n\times n}$, we write $A\succ0$ to say that $A$ is positive definite, while $A\succeq0$ means positive semi-definite. 
For $A, B\in\realR^{n\times m}$, the expression $\langle A,B\rangle$ denotes the trace of $A^\top B$. 
Given $x\in\realR^n$ and $A\in\realR^{n\times n}$, the notation $|x|^2_A$ means $x^\top Ax$. Similarly, given $B\in\realR^{m\times n}$ and $A\in\realR^{n\times n}$, the trace of $B^\top AB$ is denoted $\|B\|^2_A$.

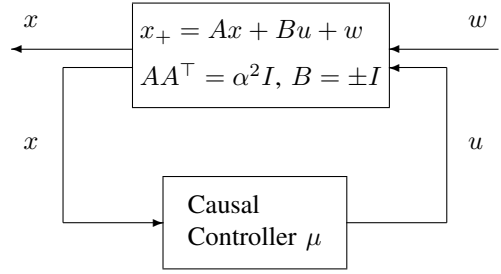
\begin{figure}
\begin{center}
\setlength{\unitlength}{.0081mm}%
\begin{picture}(9366,4600)(1318,-8461)
\put(4501,-8461){\framebox(3000,1500){\begin{tabular}{ll}Causal\\ Controller $\mu$
    \end{tabular}}}
\put(4000,-5811){\framebox(4200,1700){\begin{tabular}{l}
$x_+=Ax+Bu+w$\\[2mm]
 $AA^\top=\alpha^2I$, $B=\pm I$
\end{tabular}}}
\put(4000,-5161){\line(-1, 0){1150}}
\put(2851,-5161){\line( 0,-1){2550}}
\put(2851,-7711){\vector( 1, 0){1650}}
\put(4000,-4861){\vector(-1, 0){2000}}
\put(10001,-4861){\vector(-1, 0){1800}}
\put(7501,-7711){\line( 1, 0){1650}}
\put(9151,-7711){\line( 0, 1){2550}}
\put(9151,-5161){\vector(-1, 0){950}}
\put(2200,-6500){$x$}%
\put(9500,-6500){$u$}%
\put(2200,-4500){$x$}%
\put(9500,-4500){$w$}%
\end{picture}%
\end{center}
  \caption{We want a feedback controller that works for all system parameters within the given bounds. If $\alpha>1$, even stabilization is impossible when restricting to linear time-invariant controllers. Nonlinear adaptive controllers can do much better by estimating $(A,B)$ and use the estimate for control. The purpose of this paper is to optimize such controllers using a dynamic game formulation.
  }
\label{fig:blockdiag}
\end{figure}

\goodbreak

\section{The Solution}
\label{sec:solution}

Given $\alpha>0$, define $\gamma_*:=\alpha+\sqrt{1+\alpha^2}$ and
\begin{align}
  V^1_{A,B}(x,Z)&:=|x|^2-\gamma_*^2\left\|\begin{bmatrix}I&A&B\end{bmatrix}^\top\right\|_Z^2\label{eqn:V1}\\
  V^0_{A,B}(x,Z)&:=\frac{\gamma_*^2+1}{2}|x|^2-\frac{\gamma_*^2}{2}\sum_{i=\pm1}\left\|\begin{bmatrix}I&A&iB\end{bmatrix}^\top\right\|_Z^2\label{eqn:V0}
\end{align}

\begin{thm}
{\it Given $\alpha>0$, the value of (\ref{eqn:infsup}) is finite if and only if 
$\gamma\ge\gamma_*$. 
For $\gamma=\gamma_*$, the minimal value of (\ref{eqn:infsup}) is $(\gamma_*^2+1)|x_0|^2/2$ and a minimizing control law \eqref{eqn:mu_LQ} is 
{\small\begin{align*}
  &\mu_t(x_0,\ldots,x_t,u_0,\ldots,u_{t-1})\\
  &=\arg\min_{u_t}\max_{A,B,j,x_{t+1}}V_{A,B}^j\left(x_{\tau+1},\sum_{\tau=0}^{t}\averageE
  \begin{bmatrix}-x_{\tau+1}\\x_\tau\\u_\tau\end{bmatrix}
  \begin{bmatrix}-x_{\tau+1}\\x_\tau\\u_\tau\end{bmatrix}^\top\right)
\end{align*}
}where maximization is over $j\in\{0,1\}$, $x_{t+1}\in\realR^n$ and $A,B\in\realR^{n\times n}$ with $AA^\top=\alpha^2I$ and $B=\pm I$, while minimization is over random $u_t\in\realR^n$ with finite variance.
}
\label{thm:main}
\end{thm}
\medskip

\section{Compressing Past Data}
\label{sec:Deter_Game}

To prove Theorem~\ref{thm:main}, it is useful to first restate the original problem on a form where dynamic programming can be used. We start by observing that it is natural for the control policy to accommodate the uncertainty in $A$ and $B$ by considering historical data collected in the matrix
\begin{align*}
  Z_t&=\sum_{\tau=0}^{t-1}{\small\begin{bmatrix}-x_{\tau+1}\\x_\tau\\u_\tau\end{bmatrix}
  \begin{bmatrix}-x_{\tau+1}\\x_\tau\\u_\tau\end{bmatrix}^\top}. 
\end{align*}
This gives
$\left\|\begin{bmatrix}I\;\;\,A\;\;\,B\end{bmatrix}^\top\right\|^2_{Z_t}
  =\sum_{\tau=0}^{t-1}|w_\tau|^2
$. In fact, (\ref{eqn:infsup}) is equivalent to the following:
\smallskip

\noindent\textbf{Reformulated problem:}
\emph{Given a number $\alpha>0$ and state dimension $n$, determine
	\begin{align}
	  \!\!\!\inf_\eta\max_{A,B}\sup_{v,N}\averageE\Bigg[\sum_{t=0}^N|x_t|^2
  -\gamma^2\left\|\begin{bmatrix}I&A&B\end{bmatrix}^\top\right\|^2_{Z_{N+1}}\Bigg]
	\label{eqn:infsup2}
	\end{align}
where $AA^\top=\alpha^2I$, $B=\pm I$ and $x_t, Z_t$ are generated by
\begin{align}
  &\left\{\begin{array}{ll}
    x_{t+1}=v_t\\
    Z_{t+1}=Z_t+\averageE{\small\begin{bmatrix}-v_t\\x_t\\u_t\end{bmatrix}
    \begin{bmatrix}-v_t\\x_t\\u_t\end{bmatrix}^\top}
  \end{array}\right.
\label{eqn:xZ}
\\[3mm]
  &\;\quad u_t=\eta(x_t,Z_t).
\label{eqn:eta}
\end{align}
}The value of (\ref{eqn:infsup2}) is denoted
$V_*(x_0,Z_0)$.
\goodbreak

In this formulation, $A$ and $B$ do not appear in the dynamics, only in the penalty of the final state. As a consequence, no past states are needed in the control law (\ref{eqn:eta}), only the current state $(x_t,Z_t)$. In fact, the problem is a standard zero-sum dynamic game \cite{basar1999dynamic}, which can be addressed by dynamic programming. However, the value is the same as before:

\begin{prp}\cite{rantzer2021minimax}
  \emph{The expression (\ref{eqn:infsup}) has the value $V_*(x_0,0)$. If $\eta^*$ is optimal for (\ref{eqn:infsup2}), then $\mu^*$ defined by 
  \begin{align*}
    &\mu^*_t(x_0,\ldots,x_t,u_0,\dots,u_{t-1})\\
    &:=\arg\min_{u_t}\sup_{x_{t+1}}V_*\left(x_t,\sum_{\tau=0}^{t}\,\averageE\,{\small\begin{bmatrix}-x_{\tau+1}\\x_\tau\\u_\tau\end{bmatrix}
      \begin{bmatrix}-x_{\tau+1}\\x_\tau\\u_\tau\end{bmatrix}^\top}\right)
  \end{align*}
  is optimal for (\ref{eqn:infsup})}.
\label{prp:infsup2}
\end{prp}

\goodbreak

\begin{pf*}{Proof sketch.}
  For any fixed $N\ge0$ and $x_0\in\realR^n$, the value of (\ref{eqn:infsup}) is bounded below by the expression
  \begin{align}
    \inf_\mu\sup_{A,B,w}\averageE\sum_{t=0}^N\left(|x_t|^2-\gamma^2|w_t|^2\right),
  \label{eqn:infsup_N} 
  \end{align}
  where $x_t, w_t\in\realR^n$ and the sequences $x$ and $u$ are generated according to (\ref{eqn:plant})-(\ref{eqn:mu_LQ}). 
  The value of (\ref{eqn:infsup_N}) grows monotonically with $N$ and (\ref{eqn:infsup}) is obtained in the limit. 
  A change of variables with $v_t:=x_{t+1}$ and $Z_t$ given by (\ref{eqn:xZ}) with $Z_0=0$, shows that (\ref{eqn:infsup_N}) is equal to
  {\begin{align}
      \!\!\!\inf_\mu\sup_{A,B,v}\averageE\left[\sum_{t=0}^N|x_t|^2-\gamma^2\big\|\begin{bmatrix}I\;\;A\;\;\,B\end{bmatrix}^\top\big\|^2_{Z_{N+1}}\right]
    \label{eqn:sup_vN}
    \end{align}
    }where $x,Z,u$ are generated by (\ref{eqn:xZ}) combined with (\ref{eqn:mu_LQ}). The limit of (\ref{eqn:sup_vN}) is (\ref{eqn:infsup2}), so (\ref{eqn:infsup}) and (\ref{eqn:infsup2}) must have the same value. 
\end{pf*}

Motivated by Proposition~\ref{prp:infsup2}, we now switch attention to the reformulated problem \eqref{eqn:infsup2}-\eqref{eqn:eta}, which turns out to have an explicit soliution:

\medskip

\section{Solution to the Bellman Equation}
\label{sec:Bellman}

\begin{thm}[Explicit expression for the optimal cost]\phantom{.}\\
	\emph{
		The optimal cost of (\ref{eqn:infsup2}) is finite if and only if $\gamma\ge\gamma_*$. For $\gamma=\gamma_*$ it is given by the formula 
		\begin{align}
			V_*(x,Z)&=\max_{A,B}\max\left\{V_{A,B}^{0}(x,Z),V_{A,B}^1(x,Z)\right\}.
			\label{eqn:optimalV}
		\end{align}
		where $V_{A,B}^0$ and $V_{A,B}^1$ are defined as in \eqref{eqn:V0}-\eqref{eqn:V1}.}
	\label{thm:Bellman}
\end{thm}
\medskip

The theorem follows from Theorem~\ref{thm:upper} to be proved later. The formula will be derived by value iteration, using the following remarkable min-max result, which quantifies the optimal exploration/exploitation tradeoff:

\section{Optimal Exploration/Exploitation}

\medskip 

\begin{thm}
	\emph{ Consider $x\in\realR^n$, $g\in\realR_+$ and $Y\in\realR^{n\times(3n)}$. Let $y(A,i):=\left\langle Y,\begin{bmatrix}A&iI&iA\end{bmatrix}\right\rangle$. Then
	{\small\begin{align}
			&\min_u\max_{A,i}\max\Big\{\averageE|Ax+iu|^2+y(A,i),(g+2)|Ax|^2-g\,\averageE|u|^2\Big\}\notag\\
			&=\max_{A,i}\max\left\{y(A,i),2|\alpha x|^2\right\}
			\label{eqn:lem}
		\end{align}
	}where minimization is over random vectors $u\in\realR^n$ and maximization is over $A\in\realR^{n\times n}$ with $A^\top A=\alpha^2I$ and $i=\pm1$. Let $(\hat A,\hat i):=\argmax y(A,i)$.
	Then a minimizing argument $\bar{u}$ of \eqref{eqn:lem} is given by 
	\begin{align*}
		\begin{cases}
			\bar{u}=-\hat i\hat Ax&\text{if }y(\hat A,\hat i)\ge 2|\alpha x|^2\\
			\averageE \bar{u}=\frac{-y(\hat A,\hat i)}{2|\alpha x|^2} \hat Ax,\; \averageE|\bar{u}|^2=|\alpha x|^2&\text{otherwise.}
		\end{cases}\\[-10mm]
	\end{align*}}
	\label{thm:exploration}
\end{thm}

\begin{rmk}
	The term $y(\hat A,\hat i)$ represents infomation about the uncertain parameters $A$ and $i$ based on previous data. If $y(\hat A,\hat i)\ge 2|\alpha x|^2$, the minimizing argument is regarded as the true value and the "certainty equivalence" controller $\bar{u}=-\hat i\hat Ax$ is used. If not, a randomized component is added to learn the parameters by "exploration".
\end{rmk}

\smallskip

\begin{pf}
	Define 
	$f(x,Y)$ as the value of the left hand side in \eqref{eqn:lem} and notice that
	{\small\begin{align*}
			\max_{i}\left(\averageE|Ax+iu|^2+y(A,i)\right)
			&\ge|Ax|^2+\averageE|u|^2
			=|\alpha x|^2+\averageE|u|^2\\[-6mm]
		\end{align*}
		\begin{align*}
			f(x,Y)
			&\ge\min_{u}\max_A\max\left\{|\alpha x|^2+\averageE|u|^2,(g+2)|Ax|^2-g\,\averageE|u|^2\right\}\\
			&=\min_{u}\max\left\{|\alpha x|^2+\averageE|u|^2,(g+2)|\alpha x|^2-g\,\averageE|u|^2\right\}\\
			&=2|\alpha x|^2,
		\end{align*}
	}where the minimum is attained for $\averageE|u|^2=|\alpha x|^2$. 
	Moreover, 
	\begin{align*}
		f(x,Y)&\ge\max_{A,i}\left(\averageE|Ax+iu|^2+y(A,i)\right)\ge \max_{A,i}y(A,i).
	\end{align*}
	Combining the two lower bounds gives
	\begin{align*}
		f(x,y)&\ge\max_{A,i}\max\left\{y(A,i),2|\alpha x|^2\right\}.
	\end{align*}
	For the opposite inequality, two cases need to be considered:
	
	\medskip \goodbreak
	
	\noindent\textbf{Case~1: $y(\hat A,\hat i)\ge 2|\alpha x|^2$}\hfil\break

	\noindent Let $\bar{u}:=-\hat i\hat Ax$. Then 
	{\begin{align}
		&\max\left\{|Ax+i\bar{u}|^2+y(A,i),(g+2)|Ax|^2-g|\bar{u}|^2\right\}\notag\\
		&=\max\left\{|Ax|^2+|\bar{u}|^2+2i\bar{u}^\top Ax+y(A,i),2|\alpha x|^2\right\}\notag\\
		&=\max\left\{2|\alpha x|^2+2i\bar{u}^\top Ax+y(A,i),2|\alpha x|^2\right\}\notag\\
		&\le\max\left\{2|\alpha x|^2+2\hat i\bar{u}^\top \hat Ax+y(\hat A,\hat i),2|\alpha x|^2\right\}
		\label{ineq1}\\
		&=\max\left\{y(\hat A,\hat i),2|\alpha x|^2\right\}\notag\\
		&=\max_{A,i}\max\left\{y(A,i),2|\alpha x|^2\right\}.\notag
	\end{align}
	}
	\medskip
	
	\noindent\textbf{Case~2: $0\le y(\hat A,\hat i)\le 2|\alpha x|^2$}\hfil\break
	
	\noindent Let $\bar{u}$ be random with $\averageE \bar{u}=\frac{y(\hat A,\hat i)}{2|\alpha x|^2}\hat Ax$, $\averageE|\bar{u}|^2=|\alpha x|^2$. Then
	{\begin{align}
			&\max\left\{\averageE|Ax+i\bar{u}|^2
			+y(A,i),(g+2)|Ax|^2-g\,\averageE|\bar{u}|^2\right\}\notag\\
			&=\max\left\{|Ax|^2+|\bar{u}|^2+2i\bar{u}^\top Ax+y(A,i),2|\alpha x|^2\right\}\notag\\
			&=\max\left\{2|\alpha x|^2+2i\bar{u}^\top Ax+y(A,i),2|\alpha x|^2\right\}\notag\\
			&\le\max\left\{2|\alpha x|^2,2|\alpha x|^2\right\}
			\label{ineq2}\\
			&= 2|\alpha x|^2.\notag
		\end{align}
	}Hence, in both cases	
	\begin{align*}
		f(x,y)&\le\max_{A,i}\max\left\{y(A,i),2|\alpha x|^2\right\},
	\end{align*}
	so the proof is complete (except that the inequalities \eqref{ineq1} and \eqref{ineq2} deserve further elaboration in future versions.)
\end{pf}

\medskip

\section{Value Iteration}
\label{sec:Stoch_Game}

Let $V(x,Z)$ be a real-valued function of the state $(x,Z)$. Define the operators ${\mathcal{F}}$ and ${\mathcal{F}}_u$ according to
\begin{align*}
  &{\mathcal{F}}V(x,Z)\\
  &:=\min_u\underbrace{\max_v\left\{|x|^2
     +V\left(v,Z+\averageE{\small\begin{bmatrix}-v\\x\\u\end{bmatrix}
     \begin{bmatrix}-v\\x\\u\end{bmatrix}^\top}\right)\right\}}_{\mathcal{F}_uV(x,Z)}.
\end{align*}
(with value $+\infty$ if the supremum is unbounded for all $u$) and define the sequence $V^{(0)},V^{(1)},V^{(2)}\ldots$ according to 
\begin{align*}
	V^{(0)}(x,Z)&:=-\gamma^2\min_{A,B}\big\|\begin{bmatrix}I\;\;A\;\;\,B\end{bmatrix}^\top\big\|^2_Z
	\\
	V^{(N+1)}(x,Z)&:={\mathcal{F}}V^{(N)}(x,Z).
\end{align*}
Standard dynamic programming arguments (Proposition~\ref{prp:DynP} in the Appendix) gives $$V^{(0)}\le V^{(1)}\le\ldots\le\lim_{N\to\infty} V^{(N)}=V_*.$$
This will now be used to derive the formula for $V_*$.

\section{Lower Bound on th Optimal Cost}
\label{sec:Stoch_Game}

Define the sequence $t_0\le t_1\le t_2\ldots$ by
\begin{align}
	t_{N+1}&:=1+\frac{\gamma^2\alpha^2}{\gamma^2-t_N},\qquad\;\; t_0:=0.
	\label{eqn:tdef}
\end{align}
and let 
\begin{align*}
	V^{0,N}_{A,B}(x,Z)&:=t_N|x|^2-\frac{\gamma^2}{2}\sum_{i=\pm1}
	\left\|\begin{bmatrix}I\;\;A\;\;\,iI\end{bmatrix}^\top\right\|^2_{Z}.
\end{align*}
The limit is finite if and only if 
$(\gamma^2-t_*)(t_*-1)=\gamma^2\alpha^2$ has a real valued solution $t_*$. In other words
\begin{align*}
	4\gamma^2\alpha^2&=(\gamma^2-1)^2-
	(\gamma^2+1-2t_*)^2\\
	4\gamma^2\alpha^2&\le(\gamma^2-1)^2\\
	2\gamma \alpha&\le\gamma^2-1\\
	\alpha^2+1&\le(\gamma-\alpha)^2,
\end{align*}
which happens if and only if $\gamma\ge\gamma_*$. 

\medskip\goodbreak

\begin{thm}[Lower bound]
	\emph{
		\begin{align}
			V^{(N)}&\ge\max_{A,B}\max\left\{V_{A,B}^{0,N},V_{A,B}^1\right\}
			\label{eqn:lower}
		\end{align} 
		for all $N$.
		In particular, the value of (\ref{eqn:infsup}) is infinite if $\gamma<\gamma_*$.}
	\label{thm:lower}
\end{thm}

\begin{rmk}
	Vinnicombe \cite{vinnicombe2004examples} proved for the case $n=1$ that the value of (\ref{eqn:infsup}) is finite if
	$\gamma\ge \alpha+\sqrt{1+\alpha^2}$. He also conjectured that the reverse implication holds. Theorem~\ref{thm:lower} shows that this conjecture was correct. 
\end{rmk}

\begin{pf}
	First note that 
	\begin{align*}
		V^{(1)}&={\mathcal{F}}V^{(0)}
		=\max_{A,B}V_{A,B}^1
		=\max_{A,B}\max\left\{V_{A,B}^{0,N},V_{A,B}^1\right\},
	\end{align*}
	so the inequality \eqref{eqn:lower} holds for $N=1$. It will now be proved for all $N$ by induction.
	Assume that the inequality holds for some $N\ge 1$. To prove that it holds also for $N+1$, notice that for arbitrary (random) $u\in\realR^n$ we have
	\begin{align*}
		&\mathcal{F}_uV^{(N)}(x,Z)\\
		&\ge\mathcal{F}_u\max_{A,B}\max\left\{V_{A,B}^{0,N}(x,Z),V_{A,B}^1(x,Z)\right\}\\
		&=\max\left\{\max_{A,B}\mathcal{F}_uV_{A,B}^{0,N}(x,Z),\max_{A,B}\mathcal{F}_uV_{A,B}^1(x,Z)\right\}
	\end{align*}
	\begin{align*}
		&\mathcal{F}_uV_{A,B}^{0,N}(x,Z)-|x|^2+\gamma^2\sum_{i=\pm1}
		\left\|\begin{bmatrix}I\;\;A\;\;\,iB\end{bmatrix}^\top\right\|^2_{Z}/2\\
		&=\max_v\left\{t_N|v|^2-\gamma^2\sum_{i=\pm1}\averageE|Ax+iBu-v|^2/2\right\}\\
		&=\max_v\left\{t_N|v|^2-\gamma^2\left(|Ax-v|^2+\averageE|Bu|^2\right)\right\}\\
		&=(t_N^{-1}-\gamma^{-2})^{-1}|Ax|^2-\gamma^2\averageE|u|^2
	\end{align*}
	\begin{align*}
		&\max_{B}\mathcal{F}_uV_{A,B}^1(x,Z)-|x|^2+\gamma^2\sum_{i=\pm1}
		\left\|\begin{bmatrix}I\;\;A\;\;\,iB\end{bmatrix}^\top\right\|^2_{Z}/2\\
		&\ge\frac{1}{2}\sum_{i=\pm1}\left(\mathcal{F}_uV_{A,iB}^1(x,Z)-|x|^2+\gamma^2
		\left\|\begin{bmatrix}I\;\;A\;\;\,iB\end{bmatrix}^\top\right\|^2_{Z}\right)\\
		&=\frac{1}{2}\sum_{i=\pm1}\max_v\left(|v|^2-\gamma^2\averageE|Ax+iBu-v|^2\right)\quad\\
		&=\frac{1}{2}\sum_{i=\pm1}\frac{\averageE|Ax+iBu|^2}{1-\gamma^{-2}}\\
		&=(1-\gamma^{-2})^{-1}\left(|Ax|^2+\averageE|u|^2\right).
	\end{align*}
	Hence
	\begin{align*}
		&\min_u\mathcal{F}_uV^{(N)}(x,Z)-|x|^2+\gamma^2\sum_{i=\pm1}
		\left\|\begin{bmatrix}I\;\;A\;\;\,iB\end{bmatrix}^\top\right\|^2_{Z}/2\\
		&\ge\min_u\max\left\{\frac{|Ax|^2}{t_N^{-1}-\gamma^{-2}}-\gamma^2\averageE|u|^2,\frac{|Ax|^2+\averageE|u|^2}{1-\gamma^{-2}}\right\}\\
		&=\frac{\gamma^2|Ax|^2}{\gamma^2-t_N}=(t_{N+1}-1)|x|^2\\
		&=V_{A,B}^{0,N+1}(x,Z)-|x|^2+\gamma^2\sum_{i=\pm1}
		\left\|\begin{bmatrix}I\;\;A\;\;\,iB\end{bmatrix}^\top\right\|^2_{Z}/2,
	\end{align*}
	so $V^{(N+1)}=\min_u\mathcal{F}_uV^{(N)}\ge V_{A,B}^{0,N+1}$ and the induction argument is complete.
	For the final statement, we note that $t_N\to\infty$ when $\gamma\le\gamma_*$. Hence $\lim_{N\to\infty}V^{(N)}$ is unbounded and the same is true for the optimal cost $V^*$. 
\end{pf}

\section{Upper Bound on th Optimal Cost}

\begin{thm}[Upper bound]
	\emph{If $\gamma\ge\gamma_*$, then
		{\begin{align*}
			V^{(N)}&\le\max_{A,B}\max\left\{V_{A,B}^{0},V_{A,B}^1\right\}
		\end{align*} 
	}for all $N$.
	When $N\to\infty$, combination with Theorem~\ref{thm:lower} gives 
	\begin{align*}
		V_*&=\max_{A,B}\max\left\{V_{A,B}^{0},V_{A,B}^1\right\}.
	\end{align*}
	}
	\label{thm:upper}
\end{thm}
\medskip

\medskip

\begin{pf}
  We prove the expression for $V^{(N+1)}(x,Z)$ by induction over $N$. It holds by definition for $N=0$. Assume that it holds for some $N\ge 0$. 
  Straightforward calculations give
  \begin{align*}
    &\mathcal{F}_uV_{A,B}^1(x,Z)-|x|^2+\gamma_*^2
    \left\|\begin{bmatrix}I\;\;A\;\;\,B\end{bmatrix}^\top\right\|^2_{Z}\\
    &=\max_{v}\left\{|v|^2-\gamma_*^2\,\averageE|Ax+Bu-v|^2\right\}\\
    &=\frac{\averageE|Ax+Bu|^2}{1-\gamma_*^{-2}}\\\\
    &\mathcal{F}_uV_{A,B}^{0,N}(x,Z)-|x|^2
    +\gamma_*^2\sum_{i=\pm1}\left\|\begin{bmatrix}I\;\;A\;\;\,iB\end{bmatrix}^\top\right\|^2_{Z/2}\\
    &=\max_{v}\bigg\{t_N|v|^2-\gamma_*^2\sum_{i=\pm1}\averageE|Ax+iB u-v|^2/2\bigg\}\\
    &=\max_{v}\bigg\{t_N|v|^2-\gamma_*^2|A x-v|^2-\gamma_*^2\averageE|u|^2\bigg\}\\
    &=\frac{|A x|^2}{t_N^{-1}-\gamma_*^{-2}}-\gamma_*^2\averageE|u|^2
	\end{align*}
\begin{align*}
  &V^{(N+1)}(x,Z)-|x|^2\\
  &=\mathcal{F}V^{(N)}(x,Z)-|x|^2\\
  &\le\min_u\max_{A,B}\max\left\{\mathcal{F}_uV_{A,B}^1(x,Z),\mathcal{F}_uV_{A,B}^0(x,Z)\right\}-|x|^2\\
  &=\min_u\max_{A,B}\max\bigg\{\frac{|Ax+Bu|^2}{1-\gamma^{-2}}-\gamma^2
  \left\|\begin{bmatrix}I\;\;A\;\;\,B\end{bmatrix}^\top\right\|^2_{Z},\\
  &\qquad\frac{|A x|^2}{t_*^{-1}-\gamma^{-2}}-\gamma^2|u|^2-\gamma^2\sum_{i=\pm1}\left\|\begin{bmatrix}I\;\;A\;\;\,iB\end{bmatrix}^\top\right\|^2_{Z/2}\bigg\}\\
  &=\max_{A,B}\max\left\{V_{A,B}^1(x,Z),V_{A,B}^0(x,Z)\right\}-|x|^2
\end{align*}
Here the critical last equality is given by Theorem~\ref{thm:exploration}. This completes the induction proof.
\end{pf}

\goodbreak

\section{Example}
\begin{figure}
  \centerline{\includegraphics[width=\hsize,height=.18\vsize]{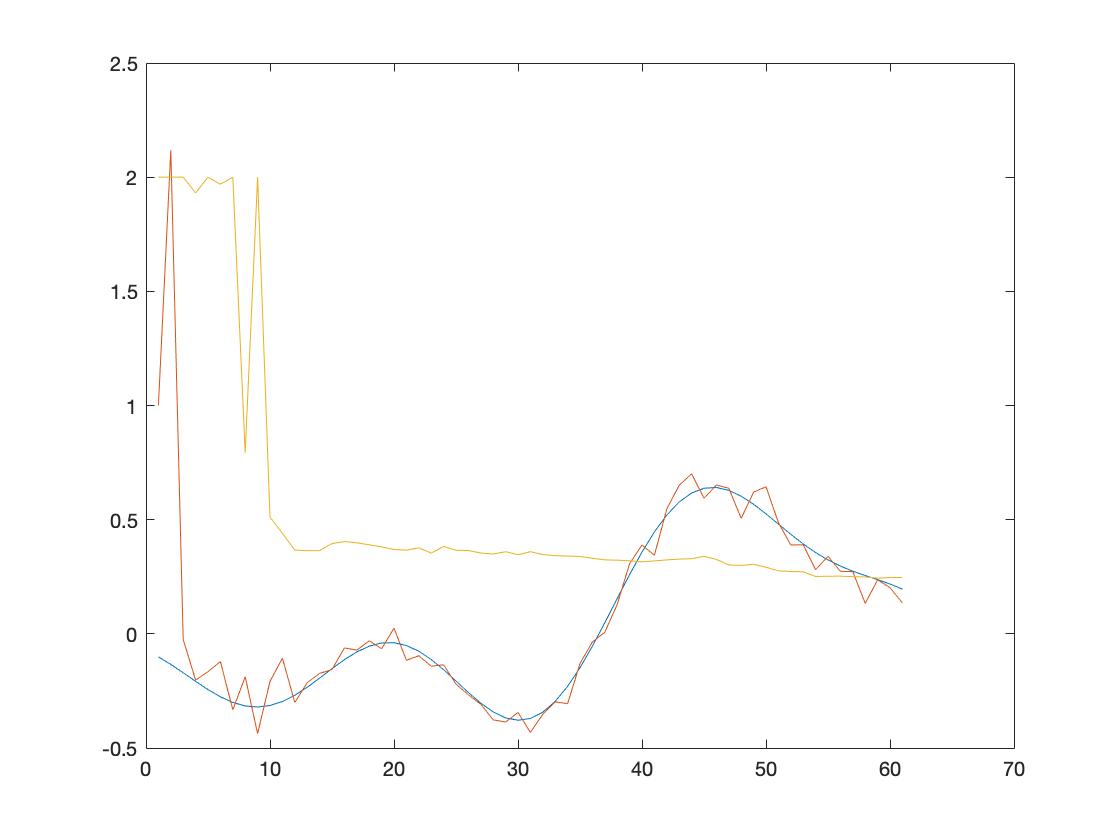}}
  \centerline{\includegraphics[width=\hsize,height=.18\vsize]{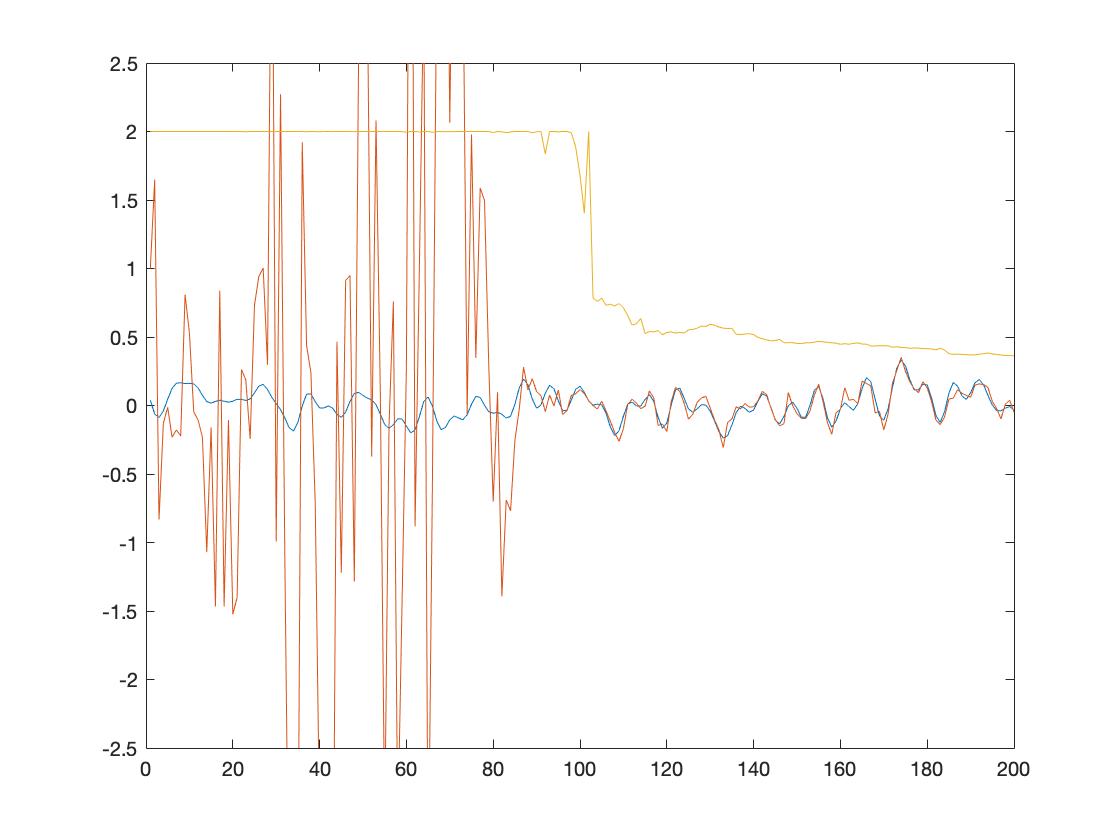}}
\caption{Two simulations of synchronizing curves are shown. The upper plot is for $n=10$, while the lower plot is for a $n=100$. The unknown $n\times n$ matrix is estimated using an optimized tradeoff between excitation and exploitation. In both figures, there is also a third curve, showing the error in the matrix estimate. The curve illustrates that the matrix convergence slows down after synchronization, since less data is then obtained. It should be noted that synchronization occurs after almost exactly $n$ time steps, when enough data has been collected to estimate all matrix directions.}
\label{fig:sync}
\end{figure}

We illustrate the theory by a simple example in $\realR^n$. A unitary matrix $A$ is randomly generated in $\realR^{n\times n}$. The we simulate the two state equations
\begin{align*}
  y_{t+1}&=Ay_t\\
  z_{t+1}&=Az_t+Bu_t+w_t
\end{align*} 
with randomly generated noise $w_t$.
To synchronize $z$ with $y$, we define $x_t=z_t-y_t$ and apply 
the policy in Theorem~\ref{thm:main}. This gives the plots shown in Figure~\ref{fig:sync} for $n=10$ and $n=100$.

\section{Conclusions}
An exact solution to the Bellman has been derived for an adaptive optimal control problem, leading to optimal tradeoff between exploration and exploitation in a dynamic setting. The problem formulation is general enough to cover state equations of arbitrary order, but several restrictions have been made. These include full state measurements, absence of control penalty and unitary $B$-matrix known except for the sign. However, the basic ideas are general, so we expect that the restrictions can be relaxed significantly in the future.


\section{Appendix: Minimax Dynamic Programming.}
  Basics facts regarding minimax dynamic programming are given below for ease reference. For more details, see for example \cite{basar1999dynamic}. 
  
  Consider the following minimax optimal control problem: Given the initial state $y_0$, consider
  \begin{align*}
    V^{(N)}(y_0)&:=\inf_\mu\sup_{v}\sum_{t=0}^{N-1}g(y_t,u_t,v_t)+g_0(y_{N})
  \end{align*}
  where supremum is taken over $v_0, v_1,\ldots, v_N$, while $y_t$ and $u_t$ are generated as    
  \begin{align*}
      y_{t+1}&=f(y_t,u_t,v_t)\\
      u_t&=\mu_t(y_0,\ldots,y_t).
  \end{align*}
  Here the control policy $\mu=(\mu_0,\mu_1,\ldots)$ generates inputs $u_t$ based on past state measurements.
  \begin{prp}
    \emph{The value function $V^{(N)}$ can be computed recursively according to $V^{(0)}(y)=g_0(y)$ and 
    \begin{align*}
      V^{(k+1)}(y)&=\min_u\max_v\left[g(y,u,v)+V^{(k)}(f(y,u,v))\right],&k&\ge1.
    \end{align*}
    Moreover, if $\max_vg(y,u,v)\ge0$ for all $(y,u)$, then 
    \begin{align}
      g_0=V^{(0)}\le V^{(1)}\le \ldots\le V^{(N)}
    \label{eqn:monotone}
    \end{align}
    and $V^{(N)}$ is upper bounded by every function $\widehat{V}$ satisfying
    \begin{align}
        g_0(y)\le \inf_u\sup_v\left[g(y,u,v)+\widehat{V}(f(y,u,v))\right]
        \le\widehat{V}(y)
    \label{eqn:Vhat}
    \end{align}
    for all $y$.}
  \label{prp:DynP}
  \end{prp}
  \medskip
  
  \begin{pf*}{Proof sketch.}
    The recursive formula can be proved by induction over $k$. By definition, it holds for $k=0$. Assume that it holds for some $k$. Then 
    {\small\begin{align*}
      &V^{(k+1)}(y_0)\\
      &=\inf_{\mu}\sup_{v}\sum_{t=0}^{k}g(y_t,u_t,v_t)+g_0(y_{k+1})\\
      &=\inf_{u_0}\sup_{v_0}\Bigg[g(y_0,u_0,v_0) +\inf_{\mu}\sup_{v}\sum_{t=1}^{k}g(y_t,u_t,v_t)+g_0(y_{k+1})\Bigg]\\
      &=\inf_{u_0}\sup_{v_0}\left[g(y_0,u_0,v_0)+V^{(k)}(f(y_0,u_0,v_0))\right],
    \end{align*}
    }where $y_{t+1}=f(y_t,u_t,v_t)$ for $0\le t\le k$, so the formula holds also for $k+1$.
    
    The assumption $\max_vg(y,u,v)\ge0$ gives $V^{(0)}\le V^{(1)}$, so the recursion gives (\ref{eqn:monotone}). 
    
    Finally, the assumption (\ref{eqn:Vhat}) gives $V^{(0)}\le\widehat{V}$. Recursion gives $V^{(N)}\le\widehat{V}$ for all $N$.
  \end{pf*}
\end{document}